\documentstyle[12pt]{article}
\newtheorem{definition}{Definition}[section]
\newtheorem{theorem}[definition]{\bf Theorem}

\newtheorem{corollary}[definition]{\bf Corollary}
\newtheorem{proposition}[definition]{\bf Proposition}
\newtheorem{remark}[definition]{\bf Remark}

\newcommand{\be}{\begin{equation}}
\newcommand{\ee}{\end{equation}}
\newcommand{\bi}{\bibitem}

\title{On $SAP$-rings}
\author{ {\small Wu Zhixiang}\\
{\small {\it Department of Mathematics}}\\
 {\small {\it Zhejiang University, Zhejiang 310027}}\\
{\small {\it P. R. China}}}

\date{}
\begin{document}

\maketitle{}

\begin{abstract}
The rings whose simple right modules are absolutely pure are
called right $SAP$-rings. We give a new characterization of right
$SAP$ rings, right $V$ rings, and von Neumann regular rings. We
also obtain a new decomposition theory of right selfinjective von
Neumann regular rings. The relationships between $SAP$-rings,
$V$-rings, and von Neumann regular rings are explored.  Some
recent results obtained by Faith are generalized and the results
of Wu-Xia are strengthened. 
\end{abstract}

\section{Introduction}

 \par\hskip 8mm Throughout this paper, the ring $R$ is an associative ring with identity, and
 all $R$-modules are unital right modules. We call a right $R$-module $M$ absolutely pure if $M$ is
 pure in every $R$-module containing $M$. Equivently, an $R$-module $M$ is absolutely pure
 if and only if it is pure in its injective hull $E(M)$. According to Stenstrom [10], an absolutely
 pure module is also called $FP$-injective module. We now call a ring $R$ right $SAP$-ring if every
 right simple $R$-module is absolutely pure. We call a ring $R$ a right $V$-ring if every
 right simple module is injective.
 Clearly, every right $V$-ring is a right $SAP$-ring. Also, it is known that a ring is von Neumann
 regular if and only if every right module is absolutely pure [10]. Thus, every von
 Neumann regular ring is a right $SAP$ ring. However, there exist SAP rings which are not von Neumann,
 because  there exist right $V$-rings which are not von Neumann regular.
 There also exist some right $SAP$ rings which are not von Neumann regular.
 Thus, the class of right $SAP$ rings is in fact a bigger class which contains
 both the class of right $V$-rings and the class of von Neumann rings as its proper
 subclasses. The class of right $SAP$ rings has been recently studied by Wu and Xia in
 [11] and by Wu and Shum in [12]. These two papers contain some basic properties of $SAP$-rings such as every $SAP$
 ring is semiprimitive; every homomorphic image of a right $SAP$ ring is still a right $SAP$
 ring;  and etc.

In this paper, we will continue this study. In Section 2, we prove
that $R$ is a right $SAP$ ring if and only if $Ann(T)=Ann(E(T))$
for any right simple module $T$ and $R/P$ are right $SAP$ rings
for all right primitive ideals $P$, where $E(T)$ is an injective
hull of $T$;  $R$ is a right $V$ ring if and only if
$Ann(T)=Ann(E(T))$ for any right simple module $T$ and $R/P$ are
right $V$ rings for all right primitive ideals $P$; $R$ is a von
Neumann regular ring if and only if $Ann(T)=Ann(E(T))$ for any
right simple module $T$ and $R/P$ are von Neumann regular rings
for all  prime ideals $P$. Many known results in [6,11,12] are
generalized ( see Corollary 2.2 and Corollary 2.4). We also obtain
a new decomposition theory of right selfinjective von Neumann
regular rings in section 2, i.e., $R=R_1\oplus R_2$, where $R_1$
is a ring of a direct product of right full linear rings, the
socle of $R_2$ is equal to zero. In Section 3, we have prove that
a ring is a right $SAP$ ring if and only if $J(N)=N\cap J(M)$ for
every submodule $N$ of any module $M$ such that $M/N$ is finitely
presented. We generalize some results in [4,6]( See Corollary 3.4
and Proposition 3.5). In the final section, we study the splitting
property of right simple modules over some $SAP$ rings. We use
$VL$ property to determine when a right $SAP$ ring is a right $V$
ring(see Corollary 4.4).

Unless stated otherwise,  all  mappings between the $R$-modules
will be \linebreak $R$-homomorphisms. We denote the submodule $N$
of a module $M$ by $N\leq M$. Also we use $M$, $E(M)$, $Soc(M)$
and $Kdim(M)$ to denote the injective hull, socle, and Krull
dimension of the module $M$ respectively. For other notations and
terminologies not given in this paper, the reader is refereed to
the texts of McConnell and Robson [9] and Stenstrom [10].

\section{$SAP$ rings and von Neumann regular rings}

\par\hskip8mm  In this section, we determine the relation among right $SAP$
rings, right $V$ rings, and von Neumann regular rings. First we
prove the following very interesting result.
\begin{theorem}\hspace{-2mm}.
(1) $R$ is a right $SAP$ ring if and only if $Ann(T)=Ann(E(T))$
for any right simple module $T$ and $R/P$ is a right $SAP$ ring
for any right primitive ideal $P$.

 (2) $R$ is a right $V$ ring if and only if $Ann(T)=Ann(E(T))$
for any right simple module $T$ and $R/P$ is a right $V$ ring for
any right primitive ideal $P$.

(3)$R$ is a von Neumann regular ring if and only if
$Ann(T)=Ann(E(T))$ for any right simple module $T$ and $R/P$ is a
von Neumann regular ring for any right prime ideal $P$.

(4) $R$ is a von Neumann regular ring if and only if $R/P$ is a
von Neumann regular ring and a left flat $R$-module for any right
prime ideal $P$.
\end{theorem}

{\it Proof.}\hskip 2mm  (1) If $R$ is a right $SAP$ ring, then
$Ann(T)=Ann(E(T))$ for any right simple module $T$ and $R/P$ is a
right $SAP$ ring for all right primitive ideal by [11, theorem 3]
and [11, proposition 5]. Conversely, assume that
$Ann(T)=Ann(E(T))$ for any right simple module $T$ and $R/P$ are
right $SAP$ rings for all right primitive ideal $P$. If $S$ is a
simple right $R$-module, then the injective hull $E(S)$ of $S$ is
also a module over $R/P$, where $P$ is the annihilator ideal of
$S$. It is obvious that $E(S)$ is an injective module over $R/P$.
For any left ideal $I$ of $R$, since $R/P$ is a right $SAP$ ring,
$SI=S(I+P)=S\cap E(S)(I+P)=S\cap E(S)I$. So $S$ is an absolutely
pure $R$-module.

(2) If $R$ is a right $V$ ring, then it is a right $SAP$ ring.
Hence $Ann(T)=Ann(E(T))$ for any right simple module $T$ and $R/P$
are right $V$ rings for all right primitive ideal by [11,Theorem
3] and [11,Theorem 6]. Conversely, $R$ is a right $SAP$ ring and
$R/P$ is a right $V$ ring for all right primitive ideal. Hence $R$
is a right $V$ ring by [11,Theorem 6].

(3) If $R$ is a von Neumann regular ring, then $Ann(T)=Ann(E(T))$
for any right simple module $T$ by [11,Theorem 3]. From [7,Theorem
1.17], we know that $R/P$ is a von Neumann regular ring for any
right prime ideal $P$. Conversely, if $Ann(T)=Ann(E(T))$ for any
right simple module $T$ and $R/P$ is a von Neumann regular ring
for any prime ideal $P$, then $R$ is a right $SAP$ ring by (1).
Hence any homomorphic image of $R$ is semiprimitive [11, Theorem
3]. This leads to $I=I^2$ for any ideal $I$. Consequently, $R$ is
a von Neumann regular by [7,Corollary 1.18].

(4) We only need to prove the sufficient part. For any right
simple module $T$ over $R$, $T$ is a right simple module over
$R/Ann(T)$. Let $E'(T)$ be an injective hull of $T$ in the
category of right $R/Ann(T)$-modules. So $Ann(E'(T))=Ann(T)=0$ in
the ring $R/Ann(T)$ by (3). Since $R/Ann(T)$ is a flat left
$R$-module, $E'(T)$ is an injective $R$-module by [7,Lemma 6.17].
Hence $E'(T)$ is an injective hull of $T$ in the category of right
$R$-modules.  Hence the annihilator ideal of $Ann(E'(T))$ in the
ring $R$ is equal to $Ann(T)$. From (3) we get $R$ is a von
Neumann regular ring. \hfill$\Box$

We hope that we can use the right primitive ideal to replace of
the prime ideal in Theorem 2.1 (3). We do not know whether a ring
$R$ is a von Neumann regular if $R/P$ is von Neumann regular for
any right primitive ideal $P$ and $Ann(T)=Ann(E(T))$ for any right
simple module over $R$. A related question  is the famous
Kaplansky's conjecture, which is wether every prime ideal of any
von Neumann regular ring is primitive. It is well-known that this
Kaplansky's  conjecture is not true for all rings ( see [4]).

If this Kaplansky's conjecture is true for a ring $R$ and $R/P$ is
right Artinian for any right primitive ideal $P$, then $R$ is a
right $SAP$ ring if and only if it is right $V$ ring if and only
if $R$ is von Neumann regular by the above Theorem. Even if the
Kaplansky's conjecture is false, we can still get the following
corollary from Theorem 2.1.

\begin{corollary}\hspace{-0.08in}. Suppose $R/P$ is either right Artinian or $PI$ for any right
primitive ideal $P$. Then the following conditions are equivalent:

(1) $R$ is a right $SAP$-ring.

(2) $R$  is a right $V$ ring.

(3)  $R$ a von Neumann regular ring.

(4) every indecomposable right $R$-module is a simple module.

\end{corollary}

{\it Proof.} Suppose $R/P$ are right Artinian for all right
primitive ideals $P$ of $R$.  Then we have $(3)\Rightarrow
(1)\Leftrightarrow (2)$ by Theorem 2.1. $(2)\Rightarrow (3)$
follows from Baccella's Theorem in [1], and $(3)\Rightarrow (4)$
follows from [8,Theorem 14]. $(4)\Rightarrow (2)$ is obvious since
$E(S)$ is indecomposable for any simple right $R$-module.

Suppose that $R/P$ are $PI$ rings for all right primitive ideals
$P$ of $R$.  Let $R$ be a right $SAP$ ring and $P$ a right
primitive ideal of $R$. If $C$ is the center of the quotient ring
$R/P$, then $C$ is a von Neumann regular domain by [11, Theorem
3]. Invoking the theorem of Kaplansky in [10,Theorem 13.3.8], we
know that $R/P$ is finitely dimensional space over $C$. Hence
$R/P$ is a semisimple Artinian ring. Thus $R$ is a von Neumann
regular right $V$ ring. Thus, $(1)\Rightarrow (2)\Rightarrow (3)$
and $(1)\Rightarrow (4)$ holds by Theorem 2.1. It is obvious that
$(3)\Rightarrow (1)$. By now we complete the proof of this
corollary. \hfill$\Box$

 There exists a ring whose every  right primitive factor is
$PI$, but $R$ is not a $PI$ ring. For example, choose a field $F$,
set $F_n=M_n(F)$, the $n\times n$ matrix ring over $F$, for
$n=1,2,\cdots,$ and let $R$ be the $F$-subalgebra of $\Pi F_n$
generated by 1 and $K=\oplus F_n$. For $n=1,2,\cdots,$ let $e_n$
be the identity element of $F_n$. Then $e_n$ is a central
idempotent in $R$, and $e_nR$ is a $PI$ ring satisfying a
polynomial identity of degree $2n+1$. Consider any primitive ideal
$P$ of $R$. If $e_n\notin P$ for some $n$, then $1-e_n\in P$. In
this case, we have $P=(1-e_n)R$ and $R/P\simeq F_n$. If all
$e_n\in P$, then $K=\oplus e_nR\leq P$. Since $R/K\simeq F$, we
see in this case $R/P\simeq F$. It is obvious that $R$ is not a
$PI$ ring.

It was proved by Kaplansky  that a commutative ring is von Neumann
regular if and only if it is a $V$ ring. This celebrated result
was extended from commutative rings to left commutative rings,
that is, a ring $R$ with $Ra\subseteq aR$ for any $a\in R$, by Wu
and Xia in [11]. In [6] Faith prove that a $PI$ ring is von
Neumann regular if and only if it is a right $V$ ring. If $R$ is a
$PI$ ring,  Wu and Shum have proved that $R$ is a von Neumann
regular ring if and only if it is a right $SAP$ ring [12, Theorem
2.4]. From the above example, we know that Corollary 2.2
generalizes these results. In addition, many restrict conditions
in [12,Theorem 2.1] can be omitted. Moreover we have the following
Corollary.

\begin{corollary}\hspace{-0.08in}. (1) Let $R$ be a ring such that
$R/P$ is finitely generated module over
its center for any primitive ideal $P$. Then $R$ is a right $SAP$
ring if and only if it is a von Neumann regular ring.

(2) Suppose $R/P$ is right Artinian for any right primitive ideal
$P$. Then $R$ is a right $SAP$ ring  if and only if every left
simple module over $R$ is flat.
\end{corollary}

{\it Proof.} (1) Since $R/P$ is a right $SAP$ ring and the center
of any right $SAP$ is a von Neumann ring [11,theorem 3], the
center $C$ of $R/P$ is a field. By the assumption $R/P$ is a
finitely dimensional over $C$. Hence $R/P$ is a von Neumann
regular ring by Theorem 2.1.

(2) Now suppose that $R/P$ is right Artinian for any right
primitive ideal $P$. If $R$ is a right $SAP$ ring, then $R$ is a
von Neumann regular ring by Corollary 2.2. So every simple left
$R$-module is flat. Now suppose that every simple left $R$-module
is flat. We prove that every right simple module is absolutely
pure in the following. Let $T$ be any simple right $R$-module and
$P=Ann(T)$. Then $(R/P)^+:=Hom(R/P,Q)$ satisfies $P(R/P)^+=0$,
where $Q$ is an injective cogenerator of the category of Abelian
groups. Thus $(R/P)^+=\sum T_i$, where $T_i$ are simple left
$R$-modules. By the assumption, $(R/P)^+$ is a flat left
$R$-module. For any finitely presented right module $M$ and any
natural number $n$ we have $0=Tor_n^R((R/P)^+,M)\simeq
Hom(Ext_R^n(M,R/P),Q)$. So $Ext_R(M,R/P)=0$ for any finitely
presented right $R$-module. Consequently, $T$ is absolutely pure.
Our proof is completed. \hfill$\Box$

In his book, Xue [13] called a ring $R$ a right quo ring if every
maximal right ideal of $R$ is an ideal of $R$. From Theorem 2.1 we
can also get the following:

\begin{corollary}\hspace{-0.08in}.  Suppose $R$ is a right quo ring. Then the following are equivalent:

(1) $Ann(T)=Ann(E(T))$ for any right simple module $T$ and any
left ideal of $R$ is an ideal of $R$;

(2) $R$  is a right $SAP$ ring;

(3) $R$ is a right $V$ ring;

(4) $R$ is a von Neumann regular and reduced ring.

(5) $R/P$ is flat left $R$-module for any maximal ideal $P$ of $R$
and any left ideal of $R$ is an ideal.
\end{corollary}

{\it Proof.} For any ring we have $(4)\Rightarrow (2)$ and
$(3)\Rightarrow (2)$. $(2)\Rightarrow (3)$ by Theorem 2.1.
$(1)\Rightarrow (2)$ by [11, Proposition 5].

Now we prove $(3)\Rightarrow (4)$.

For any nonzero element $a$, if $aR+r(a)\not =R$, then there exist
a right maximal ideal $K\supseteq aR+r(a)$. Since $R/K$ is a
simple and $R/aR$ is finitely presented, every homomorphism from
$aR$ to $R/K$ can be extend to $R$. Let $f$ be a map from $aR$ to
$R/K$ defined by $f(at)=t+K$. It is easy to prove that $f$ is a
well-defined $R$-homomorphism. This $f$ can be extend to a
homomorphism from $R$ to $R/K$. Suppose $f(1)=c+K$. Then
$1+K=f(a)=ca+K$. From this we get that $1\in K$. This is
impossible. So $aR+r(a)=R$ for any $a\in R$ and $R$ is a von
Neumann  regular and reduced ring.

Next, we prove $(4)\Rightarrow (1)$.

Suppose $R$ is a von Neumann regular and reduced ring. Then
$Ann(T)=Ann(E(T))$ for any right simple module $T$ by [11, Theorem
3]. We only need prove that every left ideal of $R$ is an ideal of
$R$ . Since $R$ is reduced, every principal left ideal is
generated by a central idempotent element. Hence any left ideal of
$R$ is an ideal of $R$.

By the above proof, we get $(4)\Rightarrow (5)$.

Finally, we prove $(5)\Rightarrow (1)$.

Let $T$ be any simple right $R$-module. Suppose $P=Ann (T)$. Then
$R/P$ is a flat left $R$-module. Let $E(T)$ be an injective hull
of $T$ in the category of right $R/P$-module. Then $E(T)$ is an
injective $R$-module. So it is also an injective hull of $T$ in
the category of right $R$-modules. Since the annihilator of $E(T)$
in $R/P$ is equal to zero, $Ann(E(T))=P$ in $R$.  Hence (1) holds.
 \hfill$\Box$

If $R$ satisfies the conditions in the Corollary 2.4 and the
maximal right quotient ring of $R$ is equal to the classical right
quotient ring of $R$, then $R$ is right and left selfinjective by
[12,Corollary 3.4] and [7,Corollary 3.9]. In [7,Chapter 9],  a
decomposition theory of right selfinjective von Neumann regular
rings has been established. We give another decomposition theory
of right selfinjective von Neumann regular rings in the following:

\begin{theorem}\hspace{-2mm}. Suppose $R$ is a right $SAP$ ring.
If $R$ is right selfinjective, then $R$ is a von Neumann regular
ring and $R=R_1\oplus R_2$, where $R_1$ is a ring of a direct
product of right full linear rings, the socle of $R_2$ is equal to
zero.
\end{theorem}

 {\it Proof.}\hskip 2mm Since $R$ is right selfinjective, $R/J$ is
a right selfinjective von Neumann regular ring. In the case that
$R$ is a right $SAP$ ring, we have $J=0$ and hence $R$ is a von
Neumann regular. In the following we prove that $R$ has the given
decomposition. More generally, we can prove that $R$ has the given
decomposition if $R$ is semiprime. Let $S$ be the right socle of
$R$. Then $R=E(S)\oplus R_2$. So there exists an idempotent, say
$f$, of $R$ such that $E(S)=fR$. Thus
$R=fRf+fR(1-f)+(1-f)Rf+(1-f)R(1-f)$. Let $L$ be any right ideal of
$R$. Then $LS\subseteq L\cap S$. For any $x\in L\cap S$, there
exists an idempotent $e$ such that $x\in Re\subseteq S$ since $S$
is a sum of idempotent left ideals. Then $x=re$ for some $r\in R$.
Therefore, $x=re\cdot e\in LS$. So $S$ is pure in $_RR$. Next we
prove that $R=fR\oplus (1-f)R$ is a direct sum of rings. First we
prove $fR(1-f)=0$. On the contrary, there is $r\in R$ such that
$fr-frf\not=0$. Recall that $S$ is pure in $_RR$. Thus,
$0=(fr-frf)S=S\cap (fr-frf)R$. As $S$ is essential in $E(S)$, so
$(fr-frf)R\cap fR=0$ and hence $fr=frf$. Consequently $fR(1-f)=0$.
From this we can prove that $(1-f)Rf=0$. In fact, if $fR(1-f)=0$,
then $(1-f)Rf$ is nilpotent ideal of $R$. Hence $(1-f)Rf=0$. By
now we have prove that $R=fR\oplus (1-f)R$ is a direct sum of
rings. Since $S$ is pure in both $_{E(S)}E(S)$ and $E(S)_{E(S)}$.
By [11,Corollary 10], $E(S)$ is isomorphic to a direct product of
right full linear rings. Let $R_1=E(S)$. It is obvious that the
socle of $R_2=(1-f)R$ is equal to zero.  By now we have completed
the proof of this theorem.\hfill$\Box$

From the proof of the above theorem we obtain the following
corollary, which generalizes a theorem of Chase and Faith
[9,Theorem 9.13].

\begin{corollary}\hspace{-2mm}. A ring $R$ is isomorphic to a direct product
of right fully linear rings if $R$ is a semiprime, right
selfinjective ring and the right socle of $R$ is an essential
right ideal of $R$. A prime right selfinjective ring is either a
right fully linear ring or a ring without socle.
\end{corollary}
 {\it Proof.}\hskip 2mm Obviously.

\section{Modules over $SAP$-rings}

\par\hskip 8mm It is well-known that a ring is a right $V$ ring if and only if
the radical of any right module is equal to zero, if and only if
any right ideal is an intersection of maximal right ideals.
Similarly, we can prove the following theorem.
\begin{theorem}\hspace{-2mm}. For any ring $R$ the following are equivalent:

(1) $R$ is a right $SAP$ ring.

(2) $J(N)=J(M)\cap N$ for any submodule N of any right module M
such that the quotient module $M/N$ is finitely presented.
\end{theorem}

{\it Proof.}\hskip 2mm $(1)\Rightarrow(2)$  Let $N$ be a submodule
of any right $R$-module $M$ and $M/N$ finitely presented. Suppose
$N_1$ is a maximal submodule of $N$. Then there exists a
homomorphism $f$ from $M$ to $N/N_1$ such that the restriction of
$f$ on $N$ is equal to the canonical projection from $N$ to
$N/N_1$. This is because that $M/N$ is finitely presented and
$N/N_1$ is $FP$-injective. Thus $kerf \cap N=N_1$. Therefore, for
any maximal submodule $N_1$ of $N$, there exists a maximal
submodule $M_1$ such that $N_1=M_1\cap N$. From this result, we
can prove that we prove $J(N)=J(M)\cap N$. In fact, let $M_1$ be a
maximal submodule of $M$. Then either $N\subseteq M_1$ or
$M=M_1+N$. In the later case, $M_1+N/M_1\simeq N/(M_1\cap N)$ is a
simple module. Thus $M_1\cap N$ is a maximal submodule of $N$.
Hence $M_1\cap N$ is a maximal submodule of $N$. So $J(M)\cap
N=(\cap M_{\alpha})\cap N=J(N)$, where $M_{\alpha}$ run through
all maximal submodules of $M$.

$(2)\Rightarrow(1)$ We need to prove that every right simple
module $T$ is $FP$-injective. For this purpose, we consider the
exact sequence $0\rightarrow N\rightarrow M\rightarrow
M/N\rightarrow 0$ with $M/N$ finitely presented. Let $g$ be a
nonzero homomorphism from $N$ to $T$. Then, clearly, $N=xR+ker g$
for any $x\in N\backslash ker g$. Fix an element $x\in N\backslash
ker g$. Since $J(T)=0$ and $M/N\simeq (M/ker g)/(N/ker g)$ is
finitely presented, there exist maximal submodules
$M_{\alpha}\supseteq kerg$ ($\alpha\in \Lambda$) of $M$ such that
$kerg=(\cap_{\alpha\in \Lambda} M_{\alpha})\cap
N=\cap_{\alpha\in\Lambda} (M_{\alpha}\cap N)$. Since $kerg$ is a
maximal submodule of $N$, there exists a maximal submodule $M_1$
of $M$ such that $ker g=M_1\cap N$. As $x\in N$, $x\notin M_1$.
This implies that $M=xR+M_1$. Consequently, $M_1\cap N=M_1\cap
(xR+N_1)=N_1+(M_1\cap xR)=N_1$. Therefore, by defining $f$ to be
zero on $M_1$ and $f=g$ on $N$, we obtain the desired extension
$g$ to $M$. \hfill$\Box$

 It is known that a ring is a right $V$ ring
if and only if the Jacobson radical of any cyclic right module is
equal to zero. Using this result and theorem 3.1, we can prove the
following result:
\begin{corollary}\hspace{-2mm}. Suppose that every simple module over the ring $R$
is finitely presented and that every submodule of any cyclic right
module has a maximal submodule. Then $R$ is a right $SAP$ ring if
and only if $R$ is a right $V$ ring. \end{corollary}

{\it Proof.}\hskip 2mm We only need to prove the sufficiency part.
Let $M$ be a right cyclic $R$-module. Suppose the radical $J(M)$
of $M$ is not equal to zero. Then there exists a maximal submodule
$N$ of $J(M)$. Since $J(M)/N$ is finitely presented,
$J(N)=J(M)\cap N=N$. Because $N$ has maximal submodules,
$N\not=J(N)$. This contradiction implies that $J(M)=0$.
Consequently, $R$ is a right $V$ ring.\hfill$\Box$

From Corollary 3.2 we get that a right Noetherian right $SAP$ ring
is a right $V$ ring. We can also give another proof of this result
as follow. Since a module is absolutely pure if and only if it is
$FP$-injective, and every absolutely pure right module over a
right Noetherian ring is injective. Hence any right Noetherian
$SAP$ ring is a right $V$ ring. We write this result as a
Corollary.

\begin{corollary}\hspace{-2mm}. Suppose that $R$ is a right Noetherian ring, then $R$
is a right $SAP$ ring if and only if $R$ is a right $V$ ring.
\end{corollary}

 Faith has given two methods to prove that a von Neumann ring is
a right $V$ ring if for any right primitive ideal $P$, $R/P$ is an
Artinian ring [6,theorem 2.1]. Using Corollary 3.3 and Theorem
2.1, we can obtain the following Corollary. Since every von
Neumann regular ring is a right $SAP$ ring, this corollary
generalizes [6,theorem 2.1].

\begin{corollary}\hspace{-0.08in}. Suppose $R$ is a  ring with right Noetherian
right primitive factor rings. Then $R$ is a right $SAP$ ring if
and only if it is a right $V$ ring.
\end{corollary}

{\it Proof.} Since a right module over a right Noetherian ring is
absolutely pure if and only if it is injective, $R/P$ is a right
$SAP$ ring if and only if it is a right $V$ ring when $R/P$ is a
right Noetherian ring. Since $R/P$ is right Noetherian for any
right primitive ideal $P$ of $R$, this corollary follows from
Theorem 2.1. \hfill$\Box$

 In
[3] Boyle and Goodearl have proved that a right $V$ ring is a
right Noetherian ring if and only if it has right Krull dimension.
We can generalize this result to right $SAP$ rings.
\begin{proposition}\hspace{-0.08in}. Suppose $R$ is a right $SAP$-ring. Then $R$ is a right
Noetherian ring if and only if it has right Krull dimension.
\end{proposition}

{\it Proof.} Every right Noetherian ring has right Krull
dimension. On the other hand, if $R$ has a right Krull dimension
and is a right $SAP$ ring, then $R$ is a right Goldie ring [9,
Proposition 6.3.5]. We claim that $R$ is a right $V$ ring in this
case. In fact, if $L$ is an essential right ideal of $R$, then by
using Corollary 3.4.7 in [9], we see that $L$ is generated by a
regular element $c\in L$, that is, $L=cR$. Let $h$ be a
homomorphism from $L$ onto a simple right module $S$. Then there
exists an element $x$ in the injective hull $E(S)$ of $S$ such
that $h(a)=xa$ for any $a\in L$. This leads to $h(c)=xc\in S\cap
E(S)c=Sc$. Thereby, there exists an element $y\in S$ such that
$xc=yc$. Now, define a homomorphism $g$ from $R$ into $S$ by
$g(a)=ya$. Then, we have $g(cb)=ycb=h(cb)$ for every element $b\in
R$. This shows that $S$ is indeed an injective $R$-module. Since
$R$ is a right $V$ ring, $R$ is a right Noetherian ring by
[3,Proposition 13]. \hfill$\Box$

In the remaind  of this section, we prove the following:
\begin{proposition} Let $M$ be a right module over a right $SAP$ ring. If $MP=M$
for all prime ideals $P$, then $M=0$.
\end{proposition}

{\it Proof.} Suppose that there is a nonzero element $x\in M$, and
choose a two-sided ideal $P$ of $R$ which is maximal with
respective to the property $x\notin xP$. We claim that $P$ is a
prime ideal of $R$. In fact, if $J$ and $K$ are two-sided ideals
of $R$ which properly contain $P$, then $xR=xJ=xK$ and so
$xR=xRK=xJK$, whence $JK$ is not contained in $P$. Thus $P$ is a
prime ideal of $R$; hence $MP=M$. Let $N$ be a maximal submodule
of $xR$ which contains $xP$. Then $xR/N$ is a simple module of
$M/N$. Using the fact that $xR/N$ is pure in $M/N$, we obtain
$0=(xR/N)P= xR/N\cap (M/N)P=xR/N$. This is a contradiction.
\hfill$\Box$

\section{The splitting property of simple modules over $SAP$-rings}
\par\hskip 8mm In this section, we want to prove that some $VL$
rings are right $V$ rings if and only if they are right $SAP$
rings.

In [2], a simple module $S_R$ is said to be self-splitting in case
$Ext_R(S,S)=0$ or, equivalently, if the category of semisimple
$S$-homogeneous right $R$-module is closed by extensions, namely,
it is a (hereditary) torsion class.

\begin{proposition}\hspace{-0.08in}. Let $R$ be a right $SAP$ ring. Suppose $S$ is a self-splitting
right simple module over $R/P$, where $P=Ann (S)$. Then $S$ is
self-splitting.
\end{proposition}

{\it Proof.} Consider any exact sequence of $R$-modules
$$0\rightarrow S\rightarrow M\rightarrow S\rightarrow 0.$$Since
$P=Ann(S)$ and $M/S\simeq S$, $MP^2=0$. As $R/P^2$ is a right
$SAP$ ring, so $P=P^2$. Hence $MP=MP^2=0$. Thus $M$ is a module
over $R/P$. So the above exact sequence is splitting.
Consequently, $S$ is a self-splitting $R$-module. \hfill$\Box$

If $R$ is a right $SAP$ ring with Artinian primitive factor, then
every right simple module is injective, thus it is self-splitting.

 In [6] a ring is called a right Camillo
ring provided that $Hom_R(E(S),E(T))=0$ for any two non-isomorphic
right simple modules $S$ and $T$. Camillo ([4]) call these rings
$H$-rings and proved that a commutative ring $R$ is an $H$-ring if
and only if $R/I$ is a local ring for all colocal ideals $I$.

\begin{theorem}\hspace{-0.08in}. Let $R$ be a right $SAP$ ring. Suppose
annihilators of any two non-isomorphic simple modules are
comaximal. Then $R$ is a right Camille ring. Moreover, suppose
$M_1\subseteq E(S_1),M_2\subseteq E(S_2)$ and right simple modules
$S_1$ and $S_2$ is not isomorphic. Then $Ext_R(M_1,M_2)=0$.
\end{theorem}

{\it Proof.} Let $S_i$ ( $i=1,2$ ) be two non-isomorphic simple
modules with annihilator ideals $P_i$. Let $E(S_i)$ be the
injective hull of $S_i$. Suppose $f$ is a nonzero homomorphism
from $E(S_1)$ to $E(S_2)$. Then $S_2P_1=S_2\cap
f(E(S_1))P_1=S_2\cap f(E(S_1)P_1)=0$. Hence $P_1\subset P_2$. This
is contradict to the assumption that $P_1$ and $P_2$ are
comaximal. Hence $f=0$ and $R$ is a right Camille ring.

Let  $0\rightarrow M_1\rightarrow M\rightarrow M_2\rightarrow 0$
be an arbitrary exact. Since $M/M_1\simeq M_2$, $P_2M\subseteq
M_1$. Consequently, $P_1P_2M=P_2P_1M=0$ and $P_1(P_1M\cap
P_2M)=P_2(P_1M\cap P_2M)=0$. Thus $M=(P_1+P_2)M=P_1M\oplus P_2M$.
As $P_2M\subseteq M_1$ and $M_1\cap P_1M=0$, so $M=M_1\oplus
P_1M$. Then $P_1M\simeq M_2$ and $Ext_R(M_1,M_2)=0$.

In [6], Faith calls a ring $R$ is right $VL$ ( for $V$-like)
provided that every subdirectly irreducible injective right module
is fieldendo. For right $SAP$ rings, we have the following.

\begin{theorem}\hspace{-0.08in}. Let $R$ be a right $SAP$ ring. Then $R$ is a right $VL$ ring if only if
 $R$ with $VL$ primitive factor
rings.
\end{theorem}

{\it Proof.}$\Rightarrow )$ Let $R_1=R/P$, where $P$ is a
primitive ideal of $R$. Suppose $M$ is a subdirectly irreducible
module over $R_1$ and $R$ is a right $VL$ ring. Then
$End_{R_1}(M)=End_R(M)$ is a division. This proves that $R/P$ is a
right $VL$ ring.

$\Leftarrow )$ If $M$ is a  subdirectly irreducible right module
over $R$, then there is a simple right $R$ module $S$ such that
$S\subseteq M\subseteq E(S)$. Set $P=Ann(S)$. Since
$Ann(S)=Ann(E(S))$, $Ann(M)=P$. Hence $M$ is a right subdirectly
irreducible module over $R/P$. So $End_{R/P}(M)=End_R(M)$ is a
division. Therefore, $R$ is a right $VL$ ring. \hfill$\Box$

\begin{corollary}\hspace{-0.08in}. Let $R$ be a ring with $VL$ right primitive factor
rings. Suppose the annihilators of two non-isomorphic simple
modules are comaximal. Then $R$ is a right $SAP$ ring if and only
if it is a right $V$ ring.
\end{corollary}

{\it Proof.} By Theorem 4.2 and Theorem 4.3, $R$ is a right $VL$
and right Camillo ring if $R$ is a right $SAP$ ring. Hence $R$ is
a right $V$ ring if it is a right $SAP$ ring by [6, Theorem 5.5].
\hfill$\Box$

\begin{remark}\hspace{-0.08in}. Suppose $R$ is a right semiartinian ring with $VL$ right primitive
factor rings. If the annihilators of two non-isomorphic simple
$R$-modules are comaximal, then $R$ is a $SAP$ ring if and only if
it is a right $V$ ring if and only if $R$ is a von Neumann regular
ring by [12, Corollary 3.11] and Corollary 4.4. \end{remark}

\begin{remark}\hspace{-0.08in}. Suppose $R/P$ is either right
Artinian or $PI$ for any right primitive
ideal $P$. Then $R$ is a ring with $VL$ right primitive factor
rings. Moreover the annihilators of two non-isomorphic simple
modules are comaximal. Thus $R$ is a right $SAP$ ring if and only
if it is a right $V$ ring.  This gives another the proof
$(3)\Rightarrow(1)\Leftrightarrow (2)$ in Corollary 2.2.
\end{remark}

\section*{ACKNOWLEDGMENT}

We would like to thank  the referee for their useful comments on
this paper. We would like to thank the CSC for the support and the
Mathematics Department of Wuppertal University for the hospitality
during the year 2003/2004.

\end{document}